 \documentclass[gmj]{degruyter-journal-a}      

\title{A Topological Approach to the Nilpotent Multipliers of a Free Product}
\headlinetitle{A Topological Approach to the Nilpotent Multipliers ...}

\lastnameone{Mashayekhy}
\firstnameone{Behrooz}
\nameshortone{B.~Mashayekhy}
\addressone{Department of Pure Mathematics,\\
Center of Excellence in Analysis on Algebraic Structures,\\
Ferdowsi University of Mashhad,\\
P.O.Box 1159-91775, Mashhad, Iran.}
\countryone{Iran}
\emailone{bmashf@um.ac.ir}

\lastnametwo{Mirebrahimi}
\firstnametwo{Hanieh}
\nameshorttwo{H.~Mirebrahimi}
\addresstwo{Department of Pure Mathematics,\\
Ferdowsi University of Mashhad,\\
P.O.Box 1159-91775, Mashhad, Iran.}
\countrytwo{Iran}
\emailtwo{h$_{-}$mirebrahimi@um.ac.ir}

\lastnamethree{Vasagh}
\firstnamethree{Zohreh}
\nameshortthree{Z.~Vasagh}
\addressthree{Department of Pure Mathematics,\\
Ferdowsi University of Mashhad,\\
P.O.Box 1159-91775, Mashhad, Iran.}
\countrythree{Iran}
\emailthree{zo$_{-}$va761@stu-mail.um.ac.ir}

\abstract{In this paper, using the topological interpretation of the Baer invariant of a group $G$, $\mathcal{V}M(G)$, with respect to an arbitrary variety $\mathcal{V}$, we extend a result of Burns and Ellis (Math. Z. 226 (1997) 405-428) on the second nilpotent multiplier of a free product of two groups to
the $c$-nilpotent multipliers, for all $c\geq 1$. In particular, we show that $M^{(c)}(G\ast H)\cong M^{(c)}(G)\oplus M^{(c)}(H)$ when $G$ and $H$ are finite groups with some conditions or when $G$ and $H$ are two perfect groups.}



\researchsupported{This research was supported by a grant from Ferdowsi University of Mashhad; (No. MP91287MSH).}

\acknowledgments{The authors would like to thank the referee for the comments and suggestions that help improve the manuscript.}

\begin{document}

\section{Introduction}

Let $G\cong F/R $ be a free presentation of a group $G$, and $\mathcal{V}$ be a variety of groups defined by a set of laws $V$. Then the Baer invariant of $G$ with respect to $\mathcal{V}$, denoted by $\mathcal{V}M(G)$, is defined to be
\[\mathcal{V}M(G)\cong\frac{R\bigcap {V(F)}}{[RV^*F]},\]
where $V(F)$ is the verbal subgroup of $F$ and
$$[RV^*F]=\big<v(f_1,...,f_{i-1},f_ir,f_{i+1},...,f_n){\big(v(f_1,...,f_n)\big)}^{-1}|$$
 $$\ \ \ \ \ \ \ \ \ \ \ \ \ \ \ \ \ \ \ \ \ r\in R,{f_i}\in F,\ 1\leq i\leq n,\ v\in V,\ n\in \mathbb{N} \big>.$$

 Note that the Baer invariant of $G$ is always abelian and independent of the presentation of $G$ (see \cite{led}). In particular, if $\mathcal{V}$ is the variety of abelian groups, then the Baer invariant of $G$ is the well-known notion, the Schur multiplier of $G$, which is isomorphic to the second homology group of $G$, $H_2(G,\Bbb{Z})$ (see \cite{kar}). If $\mathcal{V}$ is the variety of nilpotent groups of class at most $c\geq1$, then the Baer invariant of the group $G$ is called the $c$-nilpotent multiplier of $G$ which is denoted by $M^{(c)}(G)$, and will be
$$M^{(c)}(G)=\frac{R\cap \gamma_{c+1}(F)}{[R,\ _cF]},$$
where $\gamma_{c+1}(F)$ is the $(c+1)$-st term of the lower central
series of $F$ and $[R,\ _1F]=[R,F], [R,\ _cF]=[[R,\
_{c-1}F],F]$, inductively.

 Burns and Ellis in \cite{bur}, using simplicial homotopy theory introduced a topological interpretation for the $c$-nilpotent multiplier of $G$ and gave an interesting formula for the second nilpotent multiplier of the free product of two groups as follows:
\[\begin{array}{ll}M^{(2)}(G\ast H)&\cong M^{(2)}(G) \oplus M^{(2)}(H) \\ &\oplus \big(M(G)\otimes H^{ab}\big)\oplus \big(G^{ab}\otimes  M(H)\big)
\ \ \ \ \ \ \ \ (I)\\ & \oplus Tor(G^{ab},H^{ab}).\end{array}\]

Also, in \cite{fr} Franco extended the above topological interpretation to the Baer invariant of a group $G$ with respect to any variety $\mathcal{V}$. In this paper, first, we give a topological proof to show that the Baer invariant functor $\mathcal{V}M(G)$ commutes with the direct limits of a directed system of groups. Second, we intend to extend the formula $(I)$ to the $c$-nilpotent multiplier of the free product of two groups for all $c\geq 1$. In particular, we show that $M^{(c)}(G\ast H)\cong M^{(c)}(G)\oplus M^{(c)}(H)$, whenever one of the following conditions holds:\\
$(i)$ $G$ and $H$ are finite abelian groups with coprime order.\\
$(ii)$ $G$ and $H$ are finite groups with $\big(|G|,|H^{ab}|\big)=\big(|G^{ab}|,|H|\big)=1$.\\
$(iii)$ $G$ and $H$ are two finite groups with $\big(|G^{ab}|,|H^{ab}|\big)=\big(|M(G)|,|H^{ab}|\big)=\big(|G^{ab}|,|M(H)|\big)=1$.\\
$(iv)$ $G$ and $H$ are two perfect groups.

\section{Preliminaries and Notation}

 In this section, we recall some basic notations and properties of simplicial groups which will be needed in the sequel. We refer the reader to \cite{cur} or \cite{goe}, for further details.
\begin{definition}
A \emph{simplicial set} $K_{.}$ is a sequence of sets $K_0, K_1, K_2,\ldots$ together with maps $d_i:K_n\rightarrow K_{n-1}$ (faces) and $s_i:K_n\rightarrow K_{n+1}$ (degeneracies), for each $ 0\leq i\leq n$, with the following conditions:
\[ \begin{array}{lcl}
  d_j d_i  &= & \ \ \ d_{i-1} d_j \ \ \quad\text{for $ j<i$} \\
   s_j s_i & =& \ \ \ s_{i+1} s_j \ \ \quad\text{for $ j\leq i$}\\
d_j s_i&=&\begin{cases}
s_{i-1} d_j  &\text{for $ j<i$;} \\
identity &\text{for $ j=i,i+1$;} \\
s_i d_{j-1}&\text{for $ j>i+1$.} \\
\end{cases} \end{array}\]
A \emph{simplicial map} $f:K_{.}\rightarrow L_{.}$ means a sequence of functions $f_n:K_n\rightarrow L_n$, such that $f\circ d_i=d_i\circ f$, that is the following diagram commutes:
\[
\begin{CD}
K_{n+1}    @<s_i<<  {K_n}  @>{d_i}>>  K_{n-1}\\
@Vf_{n+1}VV     @Vf_{n}VV         @VVf_{n-1}V\\
K_{n+1}    @<s_i<<  {K_n}  @>d_i>>  K_{n-1}.
\end{CD}
\]
\end{definition}
Similar to topological spaces, the homotopy of two simplicial maps between simplicial sets and the homotopy groups of simplicial sets are defined.
The category of simplicial sets and topological spaces can be related by two functors as follows:
\begin{itemize}
\item The \emph{geometric realization}, $|-|$, is a functor from the category of simplicial sets to the category of CW complexes.
\item The \emph{singular simplicial}, $S_*(-)$, is a functor from the category of topological spaces to the category of simplicial sets.
\end{itemize}

A simplicial set $K_{.}$ is called a \emph{simplicial group} if each $K_i$ is a group and all faces and degeneracies are homomorphisms.

There is a basic property of simplicial groups which is due to Moore \cite{moo}, its homotopy groups $\pi_*(G_{.})$ can be obtained as the homology of a certain chain complex $(NG_{.},\partial)$.
 \begin{definition}
If $G_{.}$ is a simplicial group, then the \emph{Moore complex} $(NG_{.},\partial)$ of $G_{.}$ is the (nonabelian) chain complex defined by $(NG)_n=\cap_{i=0}^{n-1}\ker d_i$ with $\partial_n:NG_n\rightarrow~ NG_{n-1}$ which is a restriction of $d_n$.
\end{definition}
A simplicial group $G_{.}$ is said to be \emph{free} if each $G_n$ is a free group and degeneracy homomorphisms $s_i$'s send the free basis of $G_n$ into the free basis for $G_{n+1}$.
\begin{definition}
For a reduced simplicial set $K_{.}$ (i.e. $K_0=*$), let $\mathbb{G}K_{.}$ be the simplicial group defined by $(\mathbb{G}K)_n$ which is the free group generated by $K_{n+1} \backslash s_0(K_n)$, and the face and degeneracy operators are the group homomorphisms such that
\[\begin{array}{ll}d^{\mathbb{G}K}_0k&=(d_1k)(d_0k)^{-1},\\d^{\mathbb{G}K}_ik&=d_{i+1}k, \\s^{\mathbb{G}K}_ik&=s_{i+1}k,\end{array}\]
for $i>0$ and $k\in K_{n+1}$. We can consider the above notion as a functor from reduced simplicial sets to free simplicial groups which is called \emph{Kan's functor}.
\end{definition}
\begin{definition}
A free
simplicial resolution of $G$ consists of a free simplicial group $K_{.}$ with $\pi_0(K_{.}) = G$ and
$\pi_n(K_{.}) = 0$, for all $n\geq 1$.
\end{definition}
\begin{theorem}\label{fre}
 Let $K_{.}$ and $ L_{.} $ be free simplicial resolutions of $ G$, then for every $ n\geq 0 $, $ \pi_n(T(K_{.}))\cong \pi_n(T(L_{.})) $, where $ T $ is a functor in category of groups with $ T(e)=e $ (for more details see \cite{in}).
 \end{theorem}
In the following, we recall some results that will be needed in sequel.
\begin{theorem}\label{sg}
\begin{enumerate}[\upshape (i)]
\item For every simplicial group $G_{.}$, the homotopy group $\pi_n(G_{.})$ is abelian even for $n=1$ \cite[Proposition 2.2]{cur}.\label{a}
\item Every epimorphism  between simplicial groups is a fibration \cite[Lemma 3.2]{cur}.\label{f}
\item  Let $G_{.}$ be a simplicial group, then $\pi_*(G_{.})\cong H_*(NG_{.})$ \cite[Theorem 3.7]{cur}.\label{m}
\item For every simplicial set $K_{.}$, $\mathbb{G}K_{.}\simeq\Omega|K_.|$ where $\Omega X$ is the loop space of $X$ \cite[Theorem 3.16]{cur}.\label{k}
\item Let $G_{.}$ and $H_{.}$ be simplicial Abelian groups, then $H_n\big(N(G_{.}\otimes H_{.})\big)\cong H_n\big(N(G_{.})\otimes N(H_{.})\big)$ \cite[Proposition 5.6]{cur}.\label{n}
\end{enumerate}
\end{theorem}
\section{ Topological Approach to Baer Invariants}
 Let $X=K(G,1)$ be the Eilenberg-MacLane space of $G$. Then Burns and Ellis in \cite[Proposition 4.1]{bur} presented an isomorphism $ M^{(c)}(G)\cong\pi_1\big({K_{.}}/{\gamma_{c+1}(K_{.})}\big) $, where $K_{.}$ is the free simplicial group obtained from $X$ by applying Kan's functor to the reduced singular simplicial set of $X$. Burns and Ellis's interpretation for $c=~1$ is $M(G)~\cong~\pi_1\big({K_{.}}/{\gamma_{2}(K_{.})}\big)$. Also, Kan in \cite{kan} proved that $\pi_*\big({\mathbb{G}L_.}/{\gamma_2(\mathbb{G}L_.)}\big) \cong H_{*+1}(L_.)$, where $\mathbb{G}$ is the Kan's functor. Hence $H_2(G)\cong R\cap F/[R,F]=M(G)$ which is the Hopf's formula, where $G=F/R$ is a free presentation of $G$.

 Using the above notions and similar to the Burns and Ellis's interpretation we can give a topological interpretation for the Baer invariant of a group $G$ with respect to any variety $\mathcal{V}$. We recall that the following theorem was proved categorically by Franco in \cite[Theorem 7]{fr}.
\begin{theorem}\label{main}
Let $K_{.}$ be a free simplicial resolution of $G$ and $\mathcal{V}$ be a variety of groups defined by a set of laws $V$. Then the following isomorphisms hold.
\[\begin{array}{rcl} \pi_1\big(K_{.}/V(K_{.})\big)&\cong &{\mathcal {V }} M(G)\\ \pi_0{\big(K_{.}/V(K_{.})\big)}&\cong & G/V(G). \end{array}\]
\end{theorem}
\begin{proof}
Let $G\cong F/R $ be a free presentation of $G$. Then for the simplicial group $K_{.}$ obtained by applying Kan's functor to the reduced of $S_*(X)$, we have $|K_{.}|~\simeq ~\Omega X$ using Theorem \ref{sg} \eqref{k}. Therefore ${(K_{.})}_0=F$ and ${(K_{.})}_1=R \rtimes{F}$ and $d^1_0(r,f)=f$ and $d^1_1(r,f)=rf$ (see \cite[Proposition 4.1]{bur}).
 Hence ${\big(K_{.}/V(K_{.})\big)}_0=F/V(F)$, ${\big(K_{.}/V(K_{.})\big)}_1={R/[RV^*F]}\rtimes {F/V(F)}$ and $\bar{d}^1_0$ and $\bar{d}^1_1$ are induced by $d^0_1 $ and $d^1_1$, respectively. We consider the Moore chain complexes $N\big(K_{.}/V(K_{.})\big)$ and $N\big(V(K_{.})\big)$. By Theorem \ref{sg} \eqref{m} we have $\pi_0{\big(K_{.}/V(K_{.})\big)}~\cong ~G/V(G)$ and $\pi_0\big({V(K_{.})}\big)\cong V(F)/[RV^*F]$. By Theorem \ref{sg} \eqref{f} the following exact sequence of simplicial groups is a fibration.
\[0\rightarrow V(K_{.})\rightarrow K_. \rightarrow \frac{K_{.}}{V(K_{.})}\rightarrow 0.\]
Thus it induces the long exact sequence in homotopy groups as follows:
\[\cdots\rightarrow \pi_1(K_{.}) \rightarrow \pi_1\big(\frac{K_{.}}{V(K_{.})}\big)\rightarrow \pi_0\big(V(K_{.})\big)\overset {\pi_0(\subseteq)}{\rightarrow} \pi_0(K_{.}) \rightarrow \pi_0\big(\frac{K_{.}}{V(K_{.})}\big) \rightarrow 0.\]
Also $ \pi_1(K_{.})\cong \pi_1(\Omega X)\cong \pi_2(X)=0$ and similarly $\pi_0(K_{.})\cong\pi_1(X)\cong G$. Hence $\pi_1\big(K_{.}/V(K_{.})\big)\cong \ker\big(\pi_0(\subseteq)\big)\cong V(F)\cap R/[RV^*F].$
 By Theorem \ref{sg} \eqref{k} $K_{.}$ is a free simplicial resolution of $G$, therefore by Theorem \ref{fre} the result holds.
\end{proof}
 Using the above topological interpretation of Baer invariants, we intend to study the behavior of Baer invariants with direct limits with topological approach. First, we need to find the behavior of homotopy groups of simplicial groups with respect to the direct limit.
\begin{theorem}\label{di}
Let $\{^jG_{.},\varphi^j_i; {i,j}\in J\}$ be a direct system of simplicial groups $\{^jG_{.}\}$ indexed by a directed set $J$. Then
\[\pi_n\big(\displaystyle\varinjlim_{j\in J} {}^jG_{.}\big)\cong\displaystyle\varinjlim_{j\in J}\pi_n\big(^jG_{.}\big).\]
\end{theorem}
\begin{proof}
Let $^jd^k_i:{}^jG_k\rightarrow{}^jG_{k-1}$ and $^js^k_i:{}^jG_k\rightarrow{}^jG_{k+1}$ be faces and degeneracies, for $0\leq i\leq k$.
Recall that the direct limit of simplicial groups can be considered as a simplicial group as follows
\[\begin{array}{rll}(\displaystyle\varinjlim_{j\in J} {}^jG_{.})_n&=&\displaystyle\varinjlim_{j\in J}{(^jG_{.})}_n \\ d^n_i&=&\displaystyle\varinjlim_{j\in J}(^jd^n_i) \\  s^n_i&=&\displaystyle\varinjlim_{j\in J}(^js^n_i). \end{array}\]
We have the following commutative diagram
\begin{equation}\label{d}
\begin{CD}
\displaystyle\varinjlim_{j\in J}{(^jG_{.})}_{n+1}    @<\displaystyle\varinjlim_{j\in J}(^js^n_i)<<  {\displaystyle\varinjlim_{j\in J}{(^jG_{.})}_{n}}  @>\displaystyle\varinjlim_{j\in J}(^jd^n_i)>>  \displaystyle\varinjlim_{j\in J}{(^jG_{.})}_{n-1}\\
@A({}^k\theta)_{n+1}AA     @A({}^k\theta)_{n}AA   @A({}^k\theta)_{n-1}AA\\
^kG_{n+1}    @<^ks^n_i<<  {^kG_n}  @>^kd^n_i>>  ^kG_{n-1}\\
@A(\varphi^k_l)_{n+1}AA     @A(\varphi^k_l)_{n}AA   @A(\varphi^k_l)_{n-1}AA\\
^lG_{n+1}    @<^ls^n_i<<  {^lG_n}  @>^ld^n_i>>  ^lG_{n-1}.\\
\end{CD}\end{equation}
Consider the Moore chain complex $N(\varinjlim_{j\in J} {}^jG_{.})$ as follows:
\[  \begin{CD} \cdots @>{\displaystyle\varinjlim {}^jd^3_3}>> \ker\displaystyle\varinjlim_{j\in J} {}^jd^2_0 \cap \ker\displaystyle\varinjlim_{j\in J} {}^jd^2_1@>{\displaystyle\varinjlim {}^jd^2_2}>> \ker\displaystyle\varinjlim_{j\in J} {}^jd^1_0 @>{\displaystyle\varinjlim {}^jd^1_1}>> \displaystyle\varinjlim_{j\in J} {(^jG_{.})}_0. \end{CD}\]
Since direct limit of a directed system preserves exact sequences and $$\varinjlim_{j\in J} (\ker {}^jd^i_k) \cap \varinjlim_{j\in J} (\ker  {}^jd^i_{k'})=\varinjlim_{j\in J}(\ker  {}^jd^i_k \cap \ker {}^jd^i_{k'}),$$ we have the following chain complex
\[\begin{CD}  \cdots @>{\displaystyle\varinjlim {}^jd^3_3}>>\displaystyle\varinjlim_{j\in J}(\ker\ \ {}^jd^2_0 \cap \ker \ \ {}^jd^2_1)@>{\displaystyle\varinjlim {}^jd^2_2}>> \displaystyle\varinjlim_{j\in J}  \ker{}^jd^1_0 @>{\displaystyle\varinjlim {}^jd^1_1}>> \displaystyle\varinjlim_{j\in J}{(^jG_{.})}_0.\end{CD}\]
   Hence $N(\varinjlim_{j\in J}{}^jG_{.})\cong \varinjlim_{j\in J} N(^jG_{.})$ when $J$ is a directed set. Also, homology functor preserves direct limits of directed systems of simplicial groups. Therefore, using Theorem \ref{sg} ~\eqref{m}, we have
  \[\pi_n(\displaystyle\varinjlim_{j\in J} {}^jG_{.}){\cong}H_n\big(N(\displaystyle\varinjlim_{j\in J} {}^jG_{.})\big) \cong\displaystyle\varinjlim_{j\in J}  H_n \big(N(^jG_{.})\big){\cong} \displaystyle\varinjlim_{j\in J} \pi_n(^jG_{.}).\]
\end{proof}
\begin{remark}
Note that homotopy groups do not commute with direct limits of topological spaces in general and hence
Theorem \ref{di} does not hold in the category of topological spaces. To prove this, Goodwillie \cite{good} gives the following interesting example.

 Let $S^1=\{(x,y)\in \mathbb{R}^2|x^2+y^2=1\}$ be the unit circle. Let ${A_n}=\{(x,y)\in \mathbb{R}^2|x^2+y^2=1, x\leq 1-1/n\}$ be a sequence of closed arcs in $S^1$ such that $A_n$ is
in the interior of $A_{n+1}$ and such that the union $U$ of all the $A_n$ is the
complement of a point in $S^1$.
Let $X_n$ be $S^1/A_n$.
The direct limit of the diagram of circles $X_1 \rightarrow X_2 \rightarrow ...$
is $S^1/U$, a two-point space in which only one of the points is closed. Or if one
prefers to form the colimit in the category of Hausdorff or $T_1$ spaces, then the
colimit is a point. Either way, $\pi_1$ dose not commute with the direct limit.

Note that homotopy groups preserve the direct limit of filtered based spaces (for more details see \cite[page 75]{may}).
\end{remark}

Now we are in position to give a topological proof for the following theorem which was proved algebraically in \cite{mogh}.
\begin{theorem}
Let $\{G_i,\varphi_i^j ,i \in I \}$ be a directed system of groups, then
\[\displaystyle\varinjlim_{i\in I}\mathcal{V}M(G_i)\cong {\mathcal{V}M\big(\displaystyle\varinjlim_{i\in I}(G_i)\big)}.\]
\end{theorem}
\begin{proof}
Let ${K_i}_.$ be a free simplicial group corresponding to $G_i$. By \cite[Lemma 3.2]{mogh} $\varinjlim_{i\in I}{K_i}_.$ is a free simplicial group and Theorem \ref{di} implies that
 $\varinjlim_{i\in I}{K_i}_.$ is a free simplicial resolution corresponding to $\varinjlim_{i\in I}G_i$.
  Hence we have
  \[\begin{array}{ll}{\mathcal{V}M\Big(\displaystyle\varinjlim_{i\in I}(G_i)\Big)}&\cong  \pi_1\Big(\frac{\displaystyle\varinjlim {K_i}_.}{ V \big(\displaystyle\varinjlim {K_i}_.\big)}\Big)\cong \pi_1\displaystyle\varinjlim\Big(\frac{ {K_i}_.}{V\big({K_i}_.\big)}\Big)\\&\cong\displaystyle\varinjlim\pi_1\Big(\frac{{K_i}_.}{V\big({K_i}_.\big)}\Big)\cong \displaystyle\varinjlim_{i\in I}\mathcal{V}M(G_i). \end{array}\]
  \end{proof}
\section{Main Results}
In this section, by considering the variety of nilpotent groups, we intend to compute the nilpotent multipliers of the free product of two groups.

\begin{proposition}\label{g}
Let $F=K\ast L$ be the free product of two free groups $K$ and $L$ and let $\varphi:F\rightarrow K\times L$ be the natural epimorphism. Then for all $c\geq 1$, there exists the following short exact sequence
\[0\rightarrow ker \bar{\varphi}_c \rightarrow\frac{F}{\gamma_{c+1}(F)}\stackrel{\bar{\varphi}_c}{\rightarrow}\frac{K}{\gamma_{c+1}(K)}\times\frac{L}{\gamma_{c+1}(L)}\rightarrow 0,\]
where $\ker \bar{\varphi}_c\cong \frac{[K,L]^F}{[K,L,\ _{c-1}F]^F}$ which satisfies in the following exact sequence
\[0\rightarrow \frac{[K,L,\ _{c-2}F]^F}{[K,L,\ _{c-1}F]^F} \rightarrow \ker\bar{\varphi}_c\rightarrow\frac{[K,L]^F}{[K,L,\ _{c-2}F]^F}\rightarrow 0.\]
\\Moreover, we have the following isomorphism  \[\frac{[K,L,\ _{c-2}F]^F}{[K,L,\ _{c-1}F]^F}\cong\oplus\sum_{\substack{for\ some\ i+j=c}} {\underbrace{K^{ab}\otimes...\otimes K^{ab}}_{i-times}} \otimes\displaystyle{\underbrace{L^{ab}\otimes...\otimes L^{ab}}_{j-times}.}\]
\end{proposition}


\begin{proof}
Clearly the natural epimorphism $\varphi:F\rightarrow K\times L$ induces an epimorphism \[\bar{\varphi}_c: F\ /\gamma_{c+1}(F)  \rightarrow  K/\gamma_{c+1}(K) \times L /\gamma_{c+1}(L)\] given by \[ \bar{\varphi}_c\big( \omega\gamma_{c+1}(F)\big)= \big(\omega_1 \gamma_{c+1}(K), \omega_2 \gamma_{c+1}(L)\big),\] for all $c\geq~1$, where $\varphi (\omega)=\big(\omega_1, \omega_2)$.
 Therefore, we have
 \[\ker \bar{\varphi}_c\cong \frac{[K,L]^F \gamma_{c+1}(F)}{\gamma_{c+1}(F)} \cong \frac{[K,L]^F }{[K,L]^F \cap\gamma_{c+1}(F)}\cong \frac{[K,L]^F}{[K,L,\ _{c-1}F]^F}.\]
 Hence the following exact sequence exists
 \[0\rightarrow \frac{[K,L,\ _{c-2}F]^F}{[K,L,\ _{c-1}F]^F} \rightarrow \ker\bar{\varphi}_c\rightarrow\frac{[K,L]^F}{[K,L,\ _{c-2}F]^F}\rightarrow 0.\]
  Moreover, Let $K$ and $L$ be free groups on $\{x_1, \dots , x_m\}$ and $ \{x_{m+1}, \dots , x_{m+n}\}$, respectively. Then by a theorem of Hall \cite[Theorem 11.2.4]{hall} which says that $\gamma_c(F)/\gamma_{c+1}(F)$ is a free abelian group with basis of all basic commutators of weight $c$ on $\{x_1,\ \dots,\ x_{m+n}\}$, it is easy to show that $[K,L,\ _{c-2}F]^F/[K,L,\ _{c-1}F]^F$ is a free abelian group with the basis $\bar B=\{ b[K,L,\ _{c-1}F]^F | b\in B\}$, where $B=B_1-B_2-B_3$ in which $B_1$, $B_2$, $B_3$ are the set of all basic commutators of weight $c$ on $\{x_1, \dots , x_m, \dots,x_{m+n}\}$, $\{x_1, \dots , x_m\}$ and $\{x_{m+1}, \dots , x_{m+n}\}$, respectively.
 Now by the universal property of free abelian groups and tensor products we have the following isomorphism:
 \[\frac{[K,L,\ _{c-2}F]^F}{[K,L,\ _{c-1}F]^F}\cong\oplus\sum_{\substack{for\ some\ i+j=c}} {\underbrace{K^{ab}\otimes...\otimes K^{ab}}_{i-times}} \otimes\displaystyle{\underbrace{L^{ab}\otimes...\otimes L^{ab}}_{j-times}}.\]
 Note that the number of copies in the above direct sum is the number of all basic commutators subgroups of weight $c$ on $K$ and $L$.
\end{proof}

\begin{theorem}\label{main2}
Let $G,H$ be two groups with
\[G^{ab}\otimes H^{ab}=M^{(1)}(G)\otimes H^{ab}=M^{(1)}(H)\otimes G^{ab}=Tor(G^{ab},H^{ab})=0.\] Then the following isomorphism holds, for all $c\geq1$,
\[M^{(c)}(G\ast H)\cong M^{(c)}(G)\oplus M^{(c)}(H).\]
\end{theorem}
\begin{proof}
For $c=1$, by a well-known result on Schur multiplier of the free product (see \cite[Theorem 2.6.8]{kar}), we have the following isomorphism:
\[M^{(1)}(G\ast H)\cong M^{(1)}(G)\oplus M^{(1)}(H).\]
 Now we discuss in more details on the Burns and Ellis's method in \cite[Proposition 2.13]{bur}, and extend the method to any $c\geq2$.
Let $K_{.}$ and $L_{.}$ be free simplicial groups corresponding to $X=K(G,1)$ and $Y=K(H,1)$, respectively. By van Kampen theorem we have $X\vee ~Y\cong K(G\ast H,1)$ so that the free simplicial group $F_{.}$ which obtained by applying Kan's functor to the reduced singular simplicial set of $X\vee Y$ is equal to $K_{.}\ast L_{.}$. Therefore $M^{(c)}(G\ast~ H)\cong \pi_1\big({F_{.}}/{\gamma_{c+1}(F_{.})}\big)$.
 By Proposition \ref{g}, consider the following short exact sequence of simplicial groups
\[ 0\rightarrow (\ker \bar{\varphi}_c)_{.}\rightarrow \frac{F_{.}}{\gamma_{c+1}(F_{.})}\stackrel{\bar{\varphi}_c}{\rightarrow}\frac{K_{.}}{\gamma_{c+1}(K_{.})}\times\frac{L_{.}}{\gamma_{c+1}(L_{.})}\rightarrow 0,\] where $(\ker \bar{\varphi}_c)_{.}$ is a simplicial group defined by $\big((\ker \bar{\varphi}_c)_{.}\big)_n=\ker( \bar{\varphi}_c)_n$. Theorem \ref{sg} \eqref{f} implies the following long exact sequence
\[\begin{array}{ll}\cdots&\rightarrow \pi_2 \big((\ker \bar{\varphi}_c)_{.}\big) \rightarrow \pi_2\big( \frac{F_{.}}{\gamma_{c+1}(F_{.})})\stackrel{\pi_2 (\bar{\varphi}_c)}{\rightarrow}\pi_2 (\frac{K_{.}}{\gamma_{c+1}(K_{.})}\big)\oplus\pi_2 \big(\frac{L_{.}}{\gamma_{c+1}(L_{.})}\big)\\&\rightarrow \pi_1 \big((\ker \bar{\varphi}_c)_{.}\big) \rightarrow \pi_1\big( \frac{F_{.}}{\gamma_{c+1}(F_{.})}\big)\stackrel{\pi_1 (\bar{\varphi}_c)}{\rightarrow}\pi_1 \big(\frac{K_{.}}{\gamma_{c+1}(K_{.})}\big)\oplus\pi_1\big(\frac{L_{.}}{\gamma_{c+1}(L_{.})}\big)\rightarrow \cdots .\end{array}\]

Let $\alpha^K_n:\pi_n \big({F_{.}}/{\gamma_{c+1}(F_{.})}\big)\rightarrow\pi_n \big({K_{.}} /{\gamma_{c+1}(K_{.})}\big)$ and $\alpha^L_n:\pi_n \big({F_{.}}/{\gamma_{c+1}(F_{.})}\big) \\ \rightarrow~\pi_n \big({L_{.}}/{\gamma_{c+1}(L_{.})}\big)$ be homomorphisms induced by continuous maps from $X~\vee ~Y$ to $X$ and $Y$, respectively. Since $\pi_n \big({K_{.}}/{\gamma_{c+1}(K_{.})}\big)\oplus\pi_n\big({L_{.}}/{\gamma_{c+1}(L_{.})}\big)$ is a coproduct in the category of  abelian groups, there exists a unique homomorphism $$ \psi_n:\pi_n ({K_{.}}/{\gamma_{c+1}(K_{.})})\oplus\pi_n \big({L_{.}}/{\gamma_{c+1}(L_{.})}\big)\rightarrow \pi_n\big( {F_{.}}/{\gamma_{c+1}(F_{.})}\big)$$ such that $\alpha^{K_.}_n o \psi_n =p^{K_.}_n$ and $ \alpha^{L_.}_n o \psi_n =p^{L_.}_n$, where  $p^L_n$ and $p^{K_.}_n$ are projection maps. Therefore $\psi_n o \pi_n (\bar{\varphi}_c)=id$ and consequently
\[\pi_1 (\ker\bar{\varphi}_c)\oplus\pi_1 \big({K_{.}}/{\gamma_{c+1}(K_{.})}\big)\oplus\pi_1\big({L_{.}}/{\gamma_{c+1}(L_{.})}\big)\cong\pi_1\big( {F_{.}}/{\gamma_{c+1}(F_{.})}\big).\]
By Proposition \ref{g}, we have the following exact sequence of simplicial groups \[0\rightarrow \frac{[K_.,L_.,\ _{c-2}F_.]^{F_.}}{[K_.,L_.,\ _{c-1}F_.]^{F_.}} \rightarrow (\ker\bar{\varphi}_c)_.\rightarrow\frac{[K_.,L_.]^{F_.}}{[K_.,L_.,\ _{c-2}F_.]^{F_.}}\rightarrow 0.\] Theorem \ref{sg} \eqref{f} implies the following long exact sequence of homotopy groups which in low dimension takes the following form:
\[\dots \rightarrow \pi_1(\frac{[K_.,L_.,\ _{c-2}F_.]^{F_.}}{[K_.,L_.,\ _{c-1}F_.]^{F_.}}) \rightarrow \pi_1({\ker\bar{\varphi}_c})_.\rightarrow\pi_1({\ker\bar{\varphi}_{c-1}})_.\rightarrow \cdots.\] By induction on $c$, we prove that $\pi_1({\ker\bar{\varphi}_c})_.=0$. For $c=2$, Burns and Ellis in \cite[Lemma 4.2]{bur} proved that $({\ker\bar{\varphi}_2})_.\cong K^{ab}_{.}\otimes L^{ab}_{.}$. Hence
\[\begin{array}{ll}\pi_1(\ker\bar{\varphi}_2)_.&\cong H_1\big(N(K^{ab}_{.}\otimes L^{ab}_{.})\big)\\ &\cong H_1\big( N(K^{ab}_{.})\otimes N(L^{ab}_{.})\big)\\&\cong H_1\big(N(K^{ab}_{.})\big)\otimes H_0\big(N(L^{ab}_{.})\big)\oplus H_0\big(N(K^{ab}_{.})\big)\otimes H_1\big(N(L^{ab}_{.})\big)\\&\oplus Tor\Big(H_0\big(N(K^{ab}_{.})\big),H_0\big(N(L^{ab}_{.})\big)\Big)\\& \cong M^{(1)}(G)\otimes H^{ab}\oplus M^{(1)}(H)\otimes G^{ab}\oplus Tor(G^{ab},H^{ab}).\end{array}\]
Similarly, we can prove that \[\pi_0(K^{ab}_.\otimes L^{ab}_.)\cong G^{ab}\otimes H^{ab}. \]

Now let $\pi_1({\ker\bar{\varphi}_{c-1}})_.=0$. We are going to show that $\pi_1({\ker\bar{\varphi}_c})_.=0$. Since \[\frac{[K_.,L_.,\ _{c-2}F_.]^{F_.}}{[K_.,L_.,\ _{c-1}F_.]^{F_.}}\cong\oplus\sum_{for\ some\ i+j=c} {\underbrace{K^{ab}_.\otimes...\otimes K^{ab}_.}_{i-times}} \otimes\displaystyle{\underbrace{L^{ab}_.\otimes...L^{ab}_.}_{j-times}},\]  it is enough to compute $\pi_1({\underbrace{K^{ab}_.\otimes...\otimes K^{ab}_.}_{i-times}} \otimes\displaystyle{\underbrace{L^{ab}_.\otimes...\otimes L^{ab}_.}_{j-times}})$. Since $i,j\neq0$, we have
\[\begin{array}{ll} &\ \ \ \pi_1({\underbrace{K^{ab}_.\otimes...\otimes K^{ab}_.}_{i-times}} \otimes\displaystyle{\underbrace{L^{ab}_.\otimes...\otimes L^{ab}_.}_{j-times}})\\ &\cong \pi_1(K^{ab}_.\otimes L^{ab}_.)\otimes  \pi_0({\underbrace{K^{ab}_.\otimes...\otimes K^{ab}_.}_{(i-1)-times}} \otimes\displaystyle{\underbrace{L^{ab}_.\otimes...\otimes L^{ab}_.}_{(j-1)-times}})\\ &\oplus \pi_0(K^{ab}_.\otimes L^{ab}_.)\otimes \pi_1({\underbrace{K^{ab}_.\otimes...\otimes K^{ab}_.}_{(i-1)-times}} \otimes\displaystyle{\underbrace{L^{ab}_.\otimes...\otimes L^{ab}_.}_{(j-1)-times}})\\ &\oplus Tor\big(\pi_0(K^{ab}_.\otimes L^{ab}_.), \pi_0({\underbrace{K^{ab}_.\otimes...\otimes K^{ab}_.}_{(i-1)-times}} \otimes\displaystyle{\underbrace{L^{ab}_.\otimes...\otimes L^{ab}_.}_{(j-1)-times}})\big). \end{array}\]
By the hypothesis, we have $\pi_1(K^{ab}_. \otimes L^{ab}_.)=\pi_0(K^{ab}_. \otimes L^{ab}_.)=0$. Hence $\pi_1({\ker\bar{\varphi}_c})_.=0$.
\end{proof}
\begin{corollary} Let $G$ and $H$ be two groups. Then, for all $c\geq1$, we have the following isomorphism
\[M^{(c)}(G\ast H)\cong M^{(c)}(G)\oplus M^{(c)}(H),\]
 if one of the following conditions holds:\\
$(i)$ $G$ and $H$ are two abelian groups with coprime orders.\\
$(ii)$ $G$ and $H$ are two finite groups with $\big(|G|,|H^{ab}|\big)=\big(|G^{ab}|,|H|\big)=1$.\\
$(iii)$ $G$ and $H$ are two finite groups with \[\big(|G^{ab}|,|H^{ab}|\big)=\big(|M(G)|,|H^{ab}|\big)=\big(|G^{ab}|,|M(H)|\big)=1.\]
$(vi)$ $G$ and $H$ are two perfect groups.
\end{corollary}
Note that parts $(i)-(iii)$ of the above corollary are vast generalizations of a result of the first author (see \cite[Theorem 2.5]{ma}).

\end{document}